\documentclass[11pt]{amsart}
\usepackage{amsmath,latexsym,amssymb}
\usepackage{amscd}

\newcommand{\rar}{\rightarrow}
\newcommand{\lar}{\longrightarrow}

\newtheorem{Theorem}{Theorem}[section]
\newtheorem{Lemma}[Theorem]{Lemma}
\newtheorem{Corollary}[Theorem]{Corollary}
\newtheorem{Proposition}[Theorem]{Proposition}

\newtheorem{Remark}[Theorem]{Remark}
\newtheorem{Example}[Theorem]{Example}

\def\cl{\overline}

\def\demo{\noindent{\bf Proof. }}

\def\ds{\displaystyle}

\def\gr{\mbox{\rm gr}}

\def\H{\mbox{\rm H}}

\def\l{\lambda}
\def\lar{\longrightarrow}

\def\m{\mathfrak{m}}

\def\NN{\mathbb{N}}

\def\QED{\hfill$\Box$}

\def\RR{{\bf R}}
\def\rar{\rightarrow}

\def\red{\mbox{\rm red}}

\def\rees{\mbox{$\mathcal{R}$}}

\def\spec{\mbox{\rm Spec}}

\def\surj{\twoheadrightarrow}

\begin{document}

\title[Variation of Hilbert Coefficients]
{Variation of Hilbert Coefficients}

\thanks{AMS 2010 {\em Mathematics Subject Classification}.
Primary 13A30; Secondary 13B22, 13H10, 13H15.\\
The first author is partially supported by a grant from the City University of New York PSC-CUNY Research Award Program-41.
The second author is partially supported by Grant-in-Aid for Scientific Researches (C) in Japan (19540054) and by a grant from MIMS (Meiji Institute for Advanced Study of Mathematical Sciences).  The last author is partially supported by the NSF}

\author{Laura Ghezzi}
%\thanks{$^*$ Partially supported by PSC-CUNY Research Award \#60054-38 39}
\address{Department of Mathematics\\
New York City College of Technology-Cuny\\
300 Jay Street, Brooklyn, NY 11201,  U.S.A.}
\email{lghezzi@citytech.cuny.edu}

\author{Shiro Goto}
\address{Department of Mathematics, School of Science and Technology, Meiji University, 1-1-1 Higashi-mita, Tama-ku, Kawasaki 214-8571, Japan}
\email{goto@math.meiji.ac.jp}

\author{Jooyoun Hong}
\address{Department of Mathematics\\
Southern Connecticut State University\\
501 Crescent Street, New Haven, CT 06515-1533,  U.S.A.}
\email{hongj2@southernct.edu}

\author{Wolmer V. Vasconcelos}
%\thanks{The first author is partially supported by a grant from the City University of New York PSC-CUNY Research Award Program-39. The %third author is partially supported by the NSF}
\address{Department of Mathematics \\ Rutgers University \\
110 Frelinghuysen Rd, Piscataway, NJ 08854-8019, U.S.A.}
\email{vasconce@math.rutgers.edu}

\begin{abstract}
For a Noetherian local ring $(\RR, \m)$, the first two Hilbert coefficients, $e_0$ and $e_1$,  of the $I$-adic filtration of an $\m$-primary ideal $I$ are known to code for properties of $\RR$, of the blowup of $\spec(\RR)$ along $V(I)$, and  even  of their normalizations. We give estimations for these coefficients when $I$ is enlarged (in the case of $e_1$ in the same integral  closure class) for general Noetherian local rings.
\end{abstract}

\maketitle

%\tableofcontents

\section{Introduction}

\noindent
% We will discuss primarily the last issue.
Let $(\RR, \m)$ be a Noetherian local ring of dimension
$d \ge 1$, and let $I$ be an $\m$-primary ideal.
We will consider  multiplicative,
decreasing filtrations of $\RR$ ideals, \[\mathcal{A}=\{ I_n \mid
I_0=\RR,\; I_{n+1}=I I_n,\; \forall n \gg 0 \},\] integral over the $I$-adic
filtration, conveniently coded in the corresponding Rees algebra and
its associated graded ring \[{\ds \rees(\mathcal{A}) =
\bigoplus_{n\geq 0} I_n, \quad \gr_{\mathcal{A}}(\RR) =
\bigoplus_{n\geq 0} I_n/I_{n+1}.} \] Let ${\ds \cl{\rees}=\bigoplus_{n\geq 0} \cl{I^n}}$ be the integral
closure of the Rees algebra $\rees = \rees (I)$ with $I_n = I^n$ for all $n \ge 0$, which we assume to be finite over $\rees$.

 We will
consider the Hilbert coefficients $e_i(I)$ associated to $\m$-primary
ideals $I$, for $i=0,1$. These integers play important roles in the corresponding
blowup algebras. Some of these issues have a long tradition in the
context of Cohen-Macaulay local rings,
but others are of a recent vintage for general Noetherian local rings. From
the several problem areas, we highlight the following:
\begin{itemize}
\item[(i)]
 The comparison between $e_0$ and $e_1$;

\item [(ii)]$e_1$ and normalization;
\item[(iii)] The structure of $\RR$ associated to the values of $e_1$;

\item[(iv)]  Variation of $e_i$, that is how $e_i(I)$ changes
when $I$ is enlarged.
\end{itemize}

\medskip

We are concerned here with the
last item but give brief comments on the others first.

\medskip

\noindent (i) For Cohen-Macaulay rings, an  uniform bound for $e_1(I)$ first
appeared for rings of dimension $1$ in the work of D. Kirby (\cite{Kirby63}),
\[e_1(\m)\leq {{e_0(\m)}\choose{2}}.
\] Progressively, quadratic bounds of this type were developed for
arbitrary $\m$-primary ideals in all dimensions by
several authors. As a basic source, \cite{RV08} has a systematic development of these
formulas along with a comprehensive bibliography. It also treats more
general filtrations which were helpful to us.
Among the formulas which more directly influenced the authors here,
 we single out  two developed in the work of
J. Elias (\cite{Elias5, Elias8}) and M. E. Rossi and G. Valla
(\cite{RV05, RV08}).
For an $d$-dimensional Cohen-Macaulay local ring and an $\m$-primary
ideal $I$, the first bound asserts that if $I$ is
 minimally generated by $m=\nu(I)$ elements,
\[ e_1(I)\leq {e_0(I)\choose{2}}-{{m-d}\choose{2}} -\lambda(\RR/I)+1.\]
The other bound uses the $\m$-adic order of $I$, that if $I \subset
\m^s$ and  $\overline{I}\neq \overline{\m^s}$, then
\[ e_1(I)\leq {{e_0(I)-s}\choose{2}}. \] Recently, K. Hanumanthu and C. Huneke
(\cite{Huneke11}) brought a new parameter to bear on these formulas
with their proof that
\[ e_1(I)\leq {{e_0(I)-k}\choose{2}}, \]
where $k$ is the maximal  length of chains of integrally closed
ideals between $I$ and $\m$.

%To look at the set $\mathfrak{S}(I)$ of $I$-good filtrations and
%their associated graded rings, we will focus on the role of the
%Hilbert polynomial of the Hilbert function $\l(\RR/I_{n+1})$, \[
%H_{\mathcal{A}}^1(n) \stackrel{n\gg 0}{=} P_{\mathcal{A}}^1(n) = \sum_{i=0}^{d}
%(-1)^i e_i(\mathcal{A}){{n+d-i}\choose{d-i}}, \] particularly of its
%coefficients $e_0(\mathcal{A})$ and $e_1(\mathcal{A})$.

%\bigskip

%For the $\cl{I}$-adic filtration ${\ds \mathcal{N}=\{ \cl{I}^{\; n}
%\}_{n \geq 0} }$, which is an $I$-good filtration, we have \[{\ds
%\rees(\mathcal{N}) = \bigoplus_{n\geq 0} \cl{I}^{\,n} t^n .  }\] In
%this case, we denote $e_{i}(\mathcal{N})$ by
%$\cl{e}_{i}(I)=\cl{e}_{i}$.  For convenience of exposition we discuss
%mostly $I$-adic filtrations, but return occasionally to more general
%filtrations.

\medskip

%\bigskip

\noindent (ii) Since $e_1(I)\leq \overline{e}_1(I):=e_1(\overline{\rees})$, bounds with a different
character arise. A baseline is the fact that  when $\RR$ is analytically unramified, but not
necessarilly Cohen-Macaulay, one has $\overline{e}_1(I)\geq 0$
(\cite{e1bar}). An upper bound for $\overline{e}_1(I)$ (see \cite{ni1} for
other bounds) is the following.
  Let $(\RR, \mathfrak{m})$ be a
reduced Cohen-Macaulay local ring of dimension $d$,
essentially of finite type over a perfect field, and let  $I$ be an
$\mathfrak{m}$-primary ideal. Let
 $\delta$ be a regular  element of the Jacobian ideal of $\RR$. Then
 \[e_1(I)\leq \overline{e}_1(I)
\leq \frac{t}{t + 1} \bigl[ (d-1)e_0(I) + e_0((I +\delta \RR)/\delta
\RR) \bigr],\]
where $t$ is the Cohen-Macaulay type of $\RR$.
In particular, if $\RR$ is a regular local ring
 \[\overline{e}_1(I)\leq
 \frac{  (d-1)e_0(I)}{2}.\]

\medskip

\noindent (iii) When $\RR$ is not Cohen-Macaulay, the issues become less
structured since the values of $e_1(I)$ may be negative. In fact,
using the values of $e_1(I)$ for ideals generated by systems of
parameters led to the characterization of several properties
(Cohen-Macaulay, Buchsbaum, finite cohomology)
of the ring $\RR$ itself (see \cite{GhGHOPV}, \cite{GhHV}, \cite{GO},
\cite{MV}, \cite{chern}).

\medskip

\noindent (iv) We shall now outline the   main results of this note.(We refer to \cite{icbook} for basic definitions and Rees algebras theory.) Sections 2 and 3
are organized around a list of questions about the
changes that $e_0(I)$ and $e_1(I)$ undergo when $I$ varies. An important case is
\[ e_0(J), e_1(J) \lar e_0(I), e_1(I), \quad I=(J,x).
\] Clearly the optimal baseline is that of an ideal $J$ generated by
a system of parameters, but we will consider very general cases.
 As will be seen, some relationships involve
the multiplicity  $f_0(J)$  of the special fiber.
To describe one of these estimates,  let  $(\RR, \m)$ be a Noetherian local ring of
dimension $d \ge 1$, let $J$ be an $\m$-primary ideal and let $I=(J,
h_1, \ldots, h_m)$ be  integral over $J$ of reduction number
$s=\mbox{\rm red}_{J} (I)$. Then
Theorem~\ref{e1hs} asserts that
\[ e_1(I)-e_1(J) \leq \l(\RR/(J:I)){\cdot} \left[{{m+s}\choose{s}}-1\right]{\cdot} f_0(J),\]
where $f_0(J)$ is the multiplicity of the special fiber of $\rees (J)=\bigoplus_{n \ge 0}J^n$.
We add a word of warning in reading some of the formulas with terms
like $e_1(I)-e_1(J)$. When $J$ is a minimal reduction of $I$,
$e_1(J)$ is always non-positive, according to \cite{MSV10}, and
 vanishes   when
$\RR$ is Cohen--Macaulay. In fact, for unmixed local rings the
vanishing characterizes Cohen-Macaulayness (\cite{GhGHOPV}).

\medskip

In Section 3, we address the need to link the value of $\red_{J}(I)$ to other properties of $J$. This is a well-known fact when $\RR$ is a
Cohen-Macaulay ring, but we give a general formulation in
Theorem~\ref{mainresrednCM}:
 Let $(\RR,\m)$ be a Noetherian local ring of dimension $d\geq 1$ and
infinite residue field. For an $\m$-primary ideal $I$ and a minimal reduction $J$ of $I$, there exists a minimal reduction $Q$ of $I$ such that
\[ \mbox{\rm red}_{Q}(I) \leq \max \{d{\cdot} \l(\RR/J) - 2d+1, \; 0 \}. \]

\section{Upper bounds for the variations of $e_{0}(I)$ and $e_{1}(I)$}

In our calculations we make repeated use of
 the following elementary observation.

\begin{Lemma}\label{tensor}
If $(\RR, \m)$ is a Noetherian local ring and $M$ is an $\RR$--module of finite length $\l(M)$, then
\[ \l(M\otimes N)\leq \l(M){\cdot} \nu(N)\]
for every finitely generated $\RR$--module $N$, where $\nu (N)$ denotes the minimal number of generators for $N$.
\end{Lemma}

\demo Induct on $n = \l(M)$. If $n=1$, then $M\simeq \RR/\m$ and the
assertion is clear. Suppose that $n \geq 2$ and choose an $\RR$--submodule $L$ of $M$ with $\lambda(L) = 1$. By tensoring $0\rar
L \rar M \rar M/L \rar 0$ with $N$, we get the
exact sequence \[ L \otimes N \rar M\otimes N \rar
(M/L)\otimes N \rar 0.\] Since
$\l(M/L)=n-1$, the induction hypothesis shows
\[
\l( (M/L) \otimes N)
\leq (n-1)\nu(N),
\] so that
\[
\l(M \otimes N) \leq \l(L
\otimes N) + \l(  (M/L) \otimes N) \leq (1+(n-1)){\cdot}\nu(N)= \l(M)
\nu(N).
\] \QED

\begin{Theorem} Let $(\RR, \m)$ be a Noetherian local ring of dimension $d$ and let $J \subset I=(J, h)$ be
$\m$--primary ideals of $\RR$. Then \[ e_0(J) - e_0(I) \leq \l(\RR/(J:I)){\cdot} f_0(J).\]
\end{Theorem}

\demo For $n \in \NN$, consider the following filtration:
\[\begin{array}{lll}
{\ds J^{n}=M_{0} \subset M_{1}=(M_{0},\; J^{n-1}h)} \subset \cdots &\subset& {\ds M_{r-1}=(M_{r-2},\; J^{n-r+1}h^{r-1}) }\\ && \\ &\subset&   {\ds M_{r}=(M_{r-1},\; J^{n-r}h^{r})} \\ && \\  &\subset&  {\ds  \cdots \subset M_{n}=(M_{n-1}, h^n)=I^n.}
\end{array}\]
Then we obtain
\[ \l(\RR/J^n)-\l(\RR/I^n) = \l(I^n/J^n)= \l(M_{n}/M_{0}) = \sum_{r=1}^{n}  \l(M_r/M_{r-1}). \]
For each $r$, $M_r/M_{r-1}$ is generated by the image of $h^r J^{n-r}+M_{r-1}$. Consider the natural surjection
\[ \zeta: \RR/(J: h)\otimes J^{n-r} \surj M_r/M_{r-1} =(h^r J^{n-r}+M_{r-1})/M_{r-1}. \]
%%%%%%%%%%%%%%%%%%%%%%%%%%%%%%%%%%%%%%%%%%%%%%%%%%%%%%%%%%%%
%% where $\zeta(\cl{r} \otimes x)=rh^rx+M^{r-1}$.
%%%%%%%%%%%%%%%%%%%%%%%%%%%%%%%%%%%%%%%%%%%%%%%%%%%%%%%%%%%%
Using Lemma~\ref{tensor}, we have
\[{\ds \l(M_r/M_{r-1}) \leq \l( \RR/(J:h)\otimes J^{n-r}) \leq \l(\RR/(J:h)){\cdot}\nu(J^{n-r}). }\]
It follows that
\[ \l(\RR/J^n)-\l(\RR/I^n) \leq \lambda(\RR/(J:I)){\cdot}\sum_{r=0}^{n-1} \nu(J^r).\]
 The iterated Hilbert function ${\ds \sum_{r=0}^{n-1} \nu(J^r)}$ is  of polynomial type of degree $d$ with leading (binomial)
 coefficient $f_{0}(J)$. Also, for $n\gg 0$, ${\ds  \l(\RR/J^n)-\l(\RR/I^n)}$  is the difference of two polynomials of
 degree $d$  and leading (binomial) coefficients $e_0(J)$ and $e_0(I)$. Hence
\[ e_0(J) - e_0(I)  \leq \l(\RR/(J:I)){\cdot} f_0(J).\]
\QED

\begin{Theorem}\label{e1h}
 Let $(\RR, \m )$ be a Noetherian local ring of dimension $d \ge 1$ and let $J \subset I=(J,h)$ be $\m$--primary ideals. If $h$ is integral over $J$, then
\[ e_1(I) - e_1(J) \leq \mbox{\rm red}_J(I){\cdot} \l(\RR/(J:I)){\cdot} f_0(J),\]
where $\mbox{\rm red}_J(I)$ is the reduction number of $I$ with respect to $J$.
\end{Theorem}

\demo Let $s=\mbox{\rm red}_J(I)$. Then ${\ds h^{s+1}\in JI^s }$. For $n \geq s$, we obtain the following filtration:
\[\begin{array}{lll} J^{n}=M_{0} \subset M_{1}=(M_{0},\; J^{n-1}h) \subset \cdots &\subset & M_{r}=(M_{r-1},\; J^{n-r}h^{r})
\\ && \\ &\subset & \cdots \subset M_{s}=(M_{s-1},\; J^{n-s}h^s)=I^n .
\end{array} \]
Therefore
\[{\ds \l(\RR/J^n)-\l(\RR/I^n) = \l(I^{n}/J^{n}) = \sum_{r=1}^{s} \l(M_{r}/M_{r-1}) \leq \l(\RR/(J:I)) \sum_{r=1}^s \nu(J^{n-r}). }\]

Now for $n\gg 0$, $\l(\RR/J^n)-\l(\RR/I^n)$ is the difference of two polynomials of degree $d$ and with same leading (binomial)
coefficients $e_0(J)$ and $e_0(I)$, therefore it is at most a polynomial of degree $d-1$ and leading coefficient $e_1(I)-e_1(J)$.
On the other hand, for $n\gg 0$, we have
\[{\ds \l(\RR/(J:h)){\cdot} \sum_{r=1}^s \nu(J^{n-r}) \leq \l(\RR/(J:h)) {\cdot} s {\cdot} \sum_{i=0}^{d-1} (-1)^{i} f_{i}(J){{n+d-i-2}\choose{d-i-1}}, }\]
which proves that
\[ e_1(I)-e_1(J) \leq \mbox{\rm red}_J(I){\cdot} \l(\RR/(J:I)){\cdot} f_0(J). \]
\QED

\begin{Corollary}\label{e1para}
 Let $(\RR, \m )$ be a Noetherian local ring of dimension $d \geq 1$ and infinite residue field. Let $Q \subset I=(Q,h)$ be $\m$--primary ideals such that $Q$ is a minimal reduction of $I$. Then
\[ e_1(I)  \leq \mbox{\rm red}_{Q}(I) {\cdot} \l(\RR/(Q:I)) .\]
Moreover, if $\RR$ is Gorenstein, then
\[ e_1(I)  \leq \mbox{\rm red}_{Q}(I) {\cdot} (e_{0}(I) - \l(\RR/I)) .\]
\end{Corollary}

\demo The first assertion follows from $e_1(Q) \leq 0$ \cite{MSV10} and $f_{0}(Q)=1$ for every parameter ideal $Q$.
Suppose that $\RR$ is Gorenstein. Then it is enough to show that
\[\l(\RR/(Q:I)) =  e_{0}(I) - \l(\RR/I). \]
This follows from
\[  \l(\RR/(Q:I)) = \l(\RR/Q) - \l( (Q:I)/Q ) = e_{0}(Q)  -  \l( (Q:I)/Q ) = e_{0}(Q)  - \l(\RR/I) \]
because $(Q:I)/Q$ is the canonical module of $\RR/I$.
\QED

\bigskip

\begin{Example}\label{e0Ih}{\rm (\cite[Example 7.36]{icbook})
Let $k[x,y,z]$ be the polynomial ring over an infinite field $k$.
Let $\RR=k[x,y,z]_{(x,y,z)}$ and let $J$ and $I$ be $\RR$--ideals such that
\[{\ds J=(x^{a},\; y^{b},\; z^{c}) \subset (J,\; x^{\alpha} y^{\beta} z^{\gamma})=I, }\] where
${\ds \frac{\alpha}{a} + \frac{\beta}{b} + \frac{\gamma}{c} <1}$. This inequality
ensures that $h=x^{\alpha} y^{\beta} z^{\gamma} \notin \cl{J}$.
Then we have
\[e_{0}(J) - e_{0}(I) = abc - (ab \gamma +bc \alpha + ac \beta)
                      = {\ds abc \left(1- \frac{\alpha}{a} - \frac{\beta}{b} - \frac{\gamma}{c} \right).} \]
Since ${\ds (J : I) = (J: x^{\alpha} y^{\beta} z^{\gamma}) =(x^{a-\alpha},\; y^{b-\beta},\; z^{c-\gamma}) }$,
we obtain
\[\begin{array}{lll}
\l(\RR/(J:I)){{\cdot}}f_{0}(J) &= & (a-\alpha)(b-\beta)(c-\gamma) \\ && \\
 &= & abc -bc \alpha -ac \beta - ab \gamma + a \beta \gamma + b \alpha \gamma + c \alpha \beta - \alpha \beta \gamma \\ && \\
 &= & e_{0}(J) - e_{0}(I) + \alpha \beta \gamma {\ds \left( \frac{a}{\alpha}+ \frac{b}{\beta} + \frac{c}{\gamma} -1 \right)} \\ && \\
 & > & e_{0}(J) - e_{0}(I).
\end{array}
\] Let ${\ds Q =(x^{a}-z^{c},\; y^{b}-z^{c},\;x^{\alpha} y^{\beta} z^{\gamma} ) }$ and suppose that ${\ds a > 3 \alpha,\; b > 3 \beta,\;
c > 3 \gamma }$. Note that $I=(Q,\; z^c)$.  Then $Q$ is a minimal reduction of $I$ and the reduction number $\mbox{\rm red}_{Q}(I) \leq
2$.  We can estimate $e_{1}(I)$:
\[ e_{1}(I)=e_1(I)-e_1(Q) \leq 2 \l(\RR/(Q:I)) . \]
\QED
}\end{Example}

\bigskip

Now we treat a general case of Theorem~\ref{e1h}. Let $J$ be an $\m$--primary ideal and $H=(h_{1}, \ldots, h_{m})$ a set of elements integral over $J$.
Write $I= (J, H)$, where $\nu(H) = \nu(I/J)$, and consider the difference of Hilbert functions
\begin{eqnarray*}
\l(\RR/J^n)-\l(\RR/I^n) &=&\l((J,H)^n/J^n)= \l((J^n,H J^{n-1}, \ldots, H^{n-1}J, H^n)/J^n) \\
&= & \sum_{r=1}^{n}  \l(M_r/M_{r-1}),
\end{eqnarray*}
where $M_r=(J^n,H J^{n-1}, \ldots, H^{r-1} J^{n-r+1}, H^{r} J^{n-r})$.
Note that $M_r/M_{r-1}$ is generated by the image of $H^r J^{n-r}$.
More precisely, if $ I = (J, h_1, \ldots, h_m)$, then $M_r/M_{r-1}$ is
generated by batches of elements, difficult to control.
This filtration has been used by several authors when $J$ is generated by a system
of parameters. As $\nu(I^n/J^n)$ is increasing,  the method of iterating the assertion in Theorem~\ref{e1h} tends to induce a bigger upper bound for $e_1(I)-e_1(J)$  than necessary. Instead, our formulation using the filtration above wraps it differently to accommodate our
data.

\begin{Theorem}\label{e1hs} Let $(\RR, \m)$ be a Noetherian local ring of dimension $d \geq 1$, let $J$ be an $\m$-primary ideal and let $I=(J, h_1, \ldots, h_m)$ be integral over $J$ of reduction number $s=\mbox{\rm red}_{J} I$. Then
\[ e_1(I)-e_1(J) \leq \l(\RR/(J:I)){\cdot} \left[{{m+s}\choose {s}}-1\right]{\cdot} f_0(J).\]
\end{Theorem}

\demo We have already given parts of the proof. The remaining part is
to estimate the growth of the length of ${M_r/M_{r-1}} = \left[(h_1, \ldots,h_m)^r J^{n-r}+ M_{r-1} \right]/M_{r-1} $.  We note that this module
is annihilated by $J:I$ and is generated by the `monomials' in the
$h_i$ of degree $r$, with coefficients in $J^{n-r}$.  There is a
natural surjection \[ \Phi: \RR/(J:I) \otimes \RR^{b_{r}} \otimes J^{n-r}
\lar M_r/M_{r-1},\] where
%%%%%%%%%%%%%%%%%%%%%%%%%%%%%%%%%%%%%%%%%%%%%%%%%%%%%%%%%%%%%%%%%%%%%%%%%%%%%%%%
%%  $\Phi(\cl{r} \otimes e_{i} \otimes x)= r \alpha_{i}x + M_{r-1}$ and
%%%%%%%%%%%%%%%%%%%%%%%%%%%%%%%%%%%%%%%%%%%%%%%%%%%%%%%%%%%%%%%%%%%%%%%%%%%%%%%%
${\ds b_r = {{m+r-1}\choose{r}} }$.
Therefore for $n \gg 0$,
\[\begin{array}{lll} \l(\RR/J^n)-\l(\RR/I^n) &=& {\ds
\sum_{r=1}^{s} \l(M_r/M_{r-1}) } \\ && \\ &\leq & {\ds \sum_{r=1}^{s}
\l(\RR/(J:I)){\cdot} \nu(J^n){\cdot} {{m+r-1}\choose{r}} } \\ && \\ &=& {\ds \l(\RR/(J:I)){\cdot}
\nu(J^n) {\cdot}\left[  {{m+s}\choose{s}} -1 \right], }
\end{array}
\] which completes the proof.
\QED

\bigskip

In Theorem~\ref{e1hs}, if $J$ is a minimal reduction of $I$, then it is well--known that $m=\nu(I)-\nu(J)$ does not depend on $J$ because $\m I \cap J = \m J$.
%%%%%%%%%%%%%%%%%%%%%%%%%%%%%%%%%%%%%%%%%%%%%%%%%%%%%%%%%%%%%%%%%%%%%%%%%%%%%%%%%%%%%%%%%%%%%%%%
%% [Huneke and Swanson 8.3.3, 8.3.5]
%% the mapping  $\mathcal{F}(J) \rar \mathcal{F}(I)$ is one-to-one.
%%%%%%%%%%%%%%%%%%%%%%%%%%%%%%%%%%%%%%%%%%%%%%%%%%%%%%%%%%%%%%%%%%%%%%%%%%%%%%%%%%%%%%%%%%%%
Moreover, if $\RR$ is Cohen--Macaulay, then $\l(\RR/(J:I))$ does not depend on $J$ either,  because $\l(\RR/(J:I)) = e_{0}(I) - \l(\mbox{\rm H}_{m}(I))$, where $\mbox{\rm H}_{m}(I)$ is the $m$--th Koszul homology of $I$.

\begin{Proposition}\label{f0h}
Let $(\RR, \m)$ be a Noetherian local ring of dimension $d \geq 1$, let $J$ be an $\m$-primary ideal and let $I=(J, h_1, \ldots, h_m)$ be integral over $J$ of reduction number $s=\mbox{\rm red}_{J} I$. Then
\[ f_0(I) \leq \left( 1+  \l(\RR/(J:I)){\cdot} \left[{{m+s}\choose {s}}-1\right]   \right)  {\cdot} f_0(J) .\]
\end{Proposition}

\demo  By tensoring the following exact sequence with $\RR/\m$
\[ 0 \lar J^n \lar I^n \lar I^n/J^n \lar 0,\]
we obtain
\[J^n/\m J^n \lar I^n/\m I^n \lar (I^n/J^n) \otimes \RR/\m \rar 0. \] Therefore, using Lemma~\ref{tensor}, we get
\begin{eqnarray*}
\l(I^n/\m I^n) -\l(J^n/\m J^n) &\leq &\l((I^n/J^n) \otimes \RR/\m) \\ &\leq & \l(I^n/J^n) =\l(\RR/J^n)-\l(\RR/I^n). \end{eqnarray*}
This induces the
inequalities of the leading coefficients (in degree $d-1$)
\[f_0(I) - f_{0}(J) \leq e_1(I) - e_1(J). \] Using
Theorem~\ref{e1hs}, we obtain
\[ f_0(I) - f_{0}(J) \leq e_1(I) - e_1(J) \leq \l(\RR/(J:I)){\cdot} \left[{{m+s}\choose {s}}-1\right]{\cdot} f_0(J), \]
which completes the proof.
\QED

\bigskip

\begin{Remark}{\rm
Note that the formulas for the variations of $e_1(I)$ and $f_0(I)$
require that the ideal $I$ has the same integral closure as $J$.
}\end{Remark}

\bigskip

The values of the first Hilbert coefficients  are also related to the multiplicity of certain Sally modules, according to \cite[Proposition 2.8]{C09}.
Let $(\RR, \m)$ be a Noetherian local ring of dimension $d \geq 1$ with infinite residue field. Let $I$ be an $\m$--primary ideal and $Q$ a
minimal reduction of $I$. If $\dim(S_{Q}(I))=d$ and  $H^{0}_{\m}(\RR) \subset I$, then the multiplicity $s_0(Q, I)$ of the Sally
module $S_{Q}(I)$ is
\[s_0(Q,I) = e_1(I)-e_1(Q)-e_0(I)+ \l(\RR/I).  \]

\bigskip

\begin{Corollary}\label{Sally}
Let $(\RR, \m)$ be a Noetherian local ring of dimension $d \geq 1$ with infinite residue field. Let $I$ be an $\m$--primary ideal and $Q$ a
minimal reduction of $I$. Suppose that $\dim(S_{Q}(I))=d$ and that $H^{0}_{\m}(\RR) \subset I$. Then the multiplicity $s_0(Q, I)$ of the Sally
module $S_{Q}(I)$ satisfies \[s_0(Q,I) \leq -e_{0}(I) + \l(\RR/I) +
\l(\RR/(Q:I)){\cdot}\left[{{\nu(I) -d +s}\choose {s}}-1\right], \] where $s = \mathrm{red}_Q(I)$ is
the reduction number.
\end{Corollary}

\begin{Example}{\rm  Let $\RR=k[x,y]_{(x,y)}$ where $k[x,y]$ denotes the polynomial ring over an infinite field $k$. Let $\m=(x,y)$ and $I=\m^{n}=(a_1, \ldots, a_{n}, a_{n+1})$ for some $n \geq 2$.
We assume that $Q=(a_1, a_2)$ is a minimal reduction of $I$.
Let $J=(a_1, a_2, \ldots, a_{n})$. Then since $Q \subseteq J$, $I$ is integral over $J$ with $\mathrm{red}_J(I) = 1$, because $I \ne J$ and $\mathrm{red}_Q(I) = 1$.
Using ${\ds e_{1}(I)=e_{1}(\m^n) = \frac{1}{2}n(n-1)}$ and
\[{\ds e_{1}(I) - e_{1}(J) \leq \l(\RR/ (J:a_{n+1})) f_{0}(J),   }\]
we obtain
\[{ \ds e_1(J) \geq  \frac{1}{2}n(n-1) -  \l(\RR/(J:a_{n+1})) f_{0}(J) }.    \]
}\end{Example}

\bigskip

One situation that may be amenable to further analysis is when $I = J : \m$, or
more generally $I =J : \m^s$ for some values of $s$. We refer to $I$ as a
{\em socle extension} of $J$.

\begin{Remark}\label{red1}{\rm {\bf (Reduction number one)}
Let $(\RR, \m)$ be a Noetherian local ring of dimension $d \ge 1$ with infinite residue field. Let $I$ be an $\m$--primary ideal and $Q$ a minimal reduction of $I$. Suppose that $I^2 = QI$. Then by Theorem~\ref{e1hs} we get
\[ e_1(I)-e_1(Q) \leq \l(\RR/(Q:I)) {\cdot}(\nu(I)-d) \leq \l(\RR/(Q:I)){\cdot} \l(I/Q).  \]
Suppose that $\RR$ is Cohen--Macaulay. Then since $e_1(I)=e_0(I)-\l(\RR/I) = \l(I/Q)$ (\cite[2.1]{Huneke87}), it follows that
\[  e_1(I) =e_1(I)-e_1(Q) \leq  \l(\RR/(Q:I)) {\cdot}\l(I/Q) =\l(\RR/(Q:I)){\cdot} e_1(I). \]
For example, if $\RR$ is a Cohen--Macaulay local ring that is not regular and $I  = Q: \m $, then
\[  e_1(I) =e_1(I)-e_1(Q) \leq \l(\RR/(Q:I)) {\cdot}(\nu(I)-d)  \le \l(I/Q) = e_1(I), \]
which is a case when the equality in Theorem~\ref{e1hs} holds true.
}\end{Remark}

%%%%%%%%%%%%%%%%%%%%%%%%%%%%%%%%%%%%%%%%%%%%%%%%%%%%%%%%%%%%%%%%%%%%%%%%%%%%%%%%%%%%%%%%%%%%%%%%%%5
% \begin{Example}\label{socle}{\rm
% Let $(\RR, \m)$ be a Cohen-Macaulay local ring that is not regular.
% Let $Q$ be a parameter ideal and $I  = Q: \m $.
% According to \cite{CPV1} (see also \cite[Theorem 1.106]{icbook}), the ideal $I$ is integral over $Q$ and $I^2=QI$.
% %%%%%%%%%%%%%%%%%%%%%%%%%%%%%%%%%%%%%%%%%%
% % In particular $e_0(I)=e_0(L)$.
% %%%%%%%%%%%%%%%%%%%%%%%%%%%%%%%%%%%%%%%%%%%%
% Using Theorem~\ref{e1hs}, we estimate $e_1(I)$:
%\[ e_{1}(I)= e_1(I)-e_1(Q) \leq \l(\RR/(Q:I)){\cdot} \left[{{ \nu(I)-d+1}\choose {1}}-1\right]{\cdot} f_0(Q)
%=\nu(I) - d. \] As for $f_0(I)$, we use Proposition~\ref{f0h} to get
%\[ f_{0}(I) \leq e_{1}(I) - e_{1}(Q) + f_{0}(Q)  \leq \nu(I) -d +1. \]
%}\end{Example}
%
% \bigskip
%
%
% The following are  examples where the equality holds.
%
% \begin{Example}\label{2.3}{\rm
% Let $(\RR, \m)$ be an $1$--dimensional Cohen--Macaulay local ring that is not regular.
% Let $Q=(a)$ be a parameter ideal and $I=Q : \m$. Then $I^2=Q I$ and $e_1(I) = \l(I/Q)$. Moreover, using Example~\ref{socle},
% \[  e_{1}(I)= e_1(I)-e_1(Q) \leq \nu(I) - 1 =\nu(I/Q) = \l(I/Q) = e_{1}(I).  \]
% }\end{Example}
%%%%%%%%%%%%%%%%%%%%%%%%%%%%%%%%%%%%%%%%%%%%%%%%%%%%%%%%%%%%%%%%%%%%%%%%%%%%%%%%%%%%%%%%%%%%%%%%%%%%%%%%%%

\bigskip

\begin{Example}\label{2.6}{\rm
Let $\RR$ be a Cohen--Macaulay local ring of dimension $1$ with $e_0(\RR)=2$. For every
$\m$--primary ideal $I$, there exists $a \in I$ such that $I^2=aI$
(\cite[Theorem 2.5]{SV}). Hence \[ e_1(I)  = \l(\RR/(a\RR:I)){\cdot} \left[
\nu(I)-1 \right]. \] }
\end{Example}

\bigskip

\begin{Example}\label{2.4}{\rm
Let $a$ and $\ell$ be integers such that $a \geq 4$ and $\ell \geq 2$. Let $H$ be the numerical semigroup generated by $a, a\ell-1, \{ a\ell +i \}_{1 \leq i
\leq a-3}$, and put $\RR =  k[[ t^a,\; t^{a\ell-1},\; \{ t^{a\ell
+i} \}_{1 \leq i \leq a-3} ]]$ in the formal power series ring $k[[t]]$ over a field $k$. Let $I=(t^{2a\ell-a-1},\; \{ t^{3a\ell -2a-1-i} \}_{1 \leq i
\leq a-3} ) \subsetneq \RR$ and $Q = (t^{2a\ell-a-1})
\subset I$. Then $I = \omega_{\RR}$ is a canonical ideal of $\RR$ and $Q$ is a reduction of $I$. We have  \[\m I \subseteq Q \ \ \text{and} \ \ e_1(I)=\l(\RR [t^{a\ell -a-i}\mid 1 \le i \le a-3]/\RR) = a-2 =r(\RR)
\] (\cite[Lemma 2.1]{GMP}), where $r(\RR)$ is
the Cohen--Macaulay type. Hence \[a-3 = r(\RR)-1 = \nu(I/Q) = \l(I/Q)
< e_1(I) = a-2. \] Since $r(\RR) \geq 2$, the ring $\RR$
is not Gorenstein. We have $e_1(I)=\l(I/Q)+1$, so that
\medskip
\[ e_1(I) = e_0(I) -\l(\RR/I) +1. \]
\medskip
\noindent
Therefore, thanks to \cite{Sal92}, we get $I^3=QI^2$ (hence $\mathrm{red}_Q(I) = 2$) and
\[S_Q(I) \simeq B(-1)
\]
as graded $\rees(Q)$--modules, where $\rees(Q)$ denotes the Rees algebra of $Q$, $S_Q(I)$ the Sally module of $I$ with respect to $Q$, and $B = \rees(Q)/\m \rees(Q)$.
We have

\[ a-2 = e_1(I) \leq \l(\RR/(Q:I)) \left[ {{m+s}\choose {s}} -1 \right] =
{{a-1}\choose {2}} -1,  \]
since $s = 2$ and $m = a-3$. The equality $e_1(I) = \l(\RR/(Q:I)) \left[ {{m+s}\choose {s}} -1 \right]$ holds if and only if $a=4$. When this is the case, we have $H=\left<4, 4\ell-1, 4\ell+1\right>$.  }
\end{Example}

%%%%%%%%%%%%%%%%%%%%%%%%%%%%%%%%%%%%%%%%%%%%%%%%%%%%%%%%%%%%%%%%%%%%%%%%%%5
% \begin{Example}\label{2.5}{\rm
% Suppose that $m=r=1$. Then the equality holds.
% }\end{Example}
%%%%%%%%%%%%%%%%%%%%%%%%%%%%%%%%%%%%%%%%%%%%%%%%%%%%%%%%%%%%%%%%%%%%%%%%%%%

\section{The reduction number formula}

In order to make use of Theorem~\ref{e1hs}, we need information about
the reduction number of $I$ in terms related to multiplicity. Let us
recall \cite[Theorem 2.45]{icbook}:

\begin{Theorem}\label{mainresredCM}
 Let $(\RR,\m)$ be a Cohen-Macaulay local ring of dimension $d\geq 1$ and
infinite residue field. For an $\m$-primary ideal $I$, \[ \mbox{\rm
red}(I) \leq \max \left\{\frac{d {\cdot} e_{0}(I)}{o(I)} - 2d+1 \;, \; 0 \right\}\] where $o(I)$ is the $\m$-adic order of $I$.
\end{Theorem}

To establish such a result for arbitrary Noetherian local rings, we proceed
differently.
The version of the following lemma for Cohen-Macaulay rings can be
found in \cite[Chapter 3, Theorem 1.1]{Sal78}.

\begin{Lemma} \label{lengthdim1}
Let $(\RR, \m)$ be a Noetherian local ring of dimension $1$. Let $x$ be a parameter of $\RR$. Let $E$ be a finitely generated $\RR$--module and $U$ an $\RR$--submodule of $E$. Then we have the following.

\begin{enumerate}
\item[{\rm (a)}] $\nu(U)\leq \l(\RR/(x)){\cdot} \nu(E)$. Hence  $\nu(I) \leq \l(\RR/(x))$ for every ideal $I$ of $\RR$.
\item[{\rm (b)}] If $\RR$ is Cohen--Macaulay and $x$ belongs to $\m^s$, then ${\ds \nu(U) \leq \frac{\l(\RR/(x))}{s}{\cdot} \nu(E)}$.
\end{enumerate}
\end{Lemma}

\demo (a) Let $W=\H^{0}_{\m}(E)$, $E'=E/W$, and $U'=(U+W)/W$. Then $E'$ is a Cohen-Macaulay $\RR$--module of dimension $1$ and $x$ is $E'$--regular. Moreover,
\[ \l(U'/xU') = e_{0}((x), U') \leq e_{0}((x), E') =   \l(E'/xE'). \]
Consider the following two short exact sequences:
\[
\begin{CD}
0 @>>> W         @>>> E  @>>> E' @>>> 0 \hspace{0.5 in}  \\
0 @>>> U \cap W  @>>> U  @>>> U' @>>> 0 \hspace{0.5 in}  \\
\end{CD}
\]
Then we obtain
\[
\begin{array}{lll}
\nu(U) \leq \l(U/xU) & =    & \l(U'/xU') + \l( (U \cap W)/x(U \cap W)  )  \\ && \\
                     & =    & \l(U'/xU') + \l( ( 0 :_{U \cap W}  x ) )     \\ && \\
                     & \leq & \l(U'/xU') + \l( ( 0 :_{W} x ) )     \\ && \\
                     & =    & \l(U'/xU') + \l( W/xW)               \\ && \\
                     & \leq & \l(E'/xE') + \l( W/xW )    \\ && \\
                     & =    & \l(E/xE)   \\ && \\
                     & \leq & \l(\RR/(x)){\cdot} \nu(E) . \\
\end{array}
\]

%%%%%%%%%%%%%%%%%%%%%%%%%%%%%%%%%%%%%%%%%%%%%%%%%%%%%%%%%%%%%%%%%%%%%%%%%%%%%%%%%%%%%%%%%
%% \noindent {\bf Case 2 :} Let $U$ be a submodule of a free $\RR$--module $F=\RR e \oplus \RR^{n-1}$.
%% We use induction on $n$ to show that ${\ds \nu(U) \leq \l(\RR/x\RR){\cdot} n }$. If $n=1$, it is done
%% by Case I. Suppose that $n \geq 2$. Consider the following short exact sequence
%% \[ 0 \lar \RR e \lar F=\RR e \oplus \RR^{n-1} \stackrel{\pi}{\lar} \RR^{n-1} \lar 0, \] which is split.
%% Let $U_0=U \cap \RR e \subset \RR e$ and $U'=\pi(U) \subset \RR^{n-1} $. By
%% induction hypothesis, we get
%% \[ \nu(U) \leq \nu(U_{0}) + \nu(U') \leq \l(\RR/x\RR) + \l(\RR/x\RR) {\cdot} (n-1) = \l(\RR/x\RR){\cdot} n.  \]
%%
%%
%% \medskip
%%
%% \noindent {\bf Case 3 :} Let $U$ be a submodule of $E$ with $\nu(E)=n$. Then there is a natural surjection :
%% \[  \pi : \RR^n \surj E \lar 0.   \] Then $\pi^{-1}(U)$ is a
%% submodule of $\RR^n$ so that $\nu(\pi^{-1}(U)) \leq \l(\RR/x\RR){\cdot} n$  by Case II.
%% Therefore we obtain
%% \[ \nu(U) \leq \nu(\pi^{-1}(U)) \leq \l(\RR/x\RR) {\cdot} n. \]
%%%%%%%%%%%%%%%%%%%%%%%%%%%%%%%%%%%%%%%%%%%%%%%%%%%%%%%%%%%%%%%%%%%%%%%%%%%%%%%%%%%%%%%%%%%%%%%%%%%%%%%%%%%

\medskip

\noindent (b) We may assume that the field $\RR/\m$ is infinite. Let $y \RR$ be a minimal reduction of $\m$. Then since $x \in \m^s \subseteq \overline{y^s\RR}$, we get
$$\l(\RR/(x)) = e_0(x\RR) \ge e_0(\overline{y^s\RR}) = e_0(y^s\RR)= s{{\cdot}}e_0(y\RR) = s{{\cdot}}\l(\RR/(y)).$$

Hence $\l(\RR/(y)) \le \frac{\l(\RR/(x))}{s},$ so that
\[ \nu(U) \leq \l(\RR/(y)) {\cdot} \nu(E) \leq \frac{\l(\RR/(x))}{s} {\cdot} \nu(E).\]
\QED

%%%%%%%%%%%%%%%%%%%%%%%%%%%%%%%%%%%%%%%%%%%%%%%%%%%%%%%%%%%%%%%%%%%%%%%%%%%%%%%%%%%%%%%%%%%%%%%%%%%%%%%%%%%%%%%
%%% {\bf Comment:} To be useful we must have an idea of the reduction number of $L$ in terms of the data. See the next frame.
%%%%%%%%%%%%%%%%%%%%%%%%%%%%%%%%%%%%%%%%%%%%%%%%%%%%%%%%%%%%%%%%%%%%%%%%%%%%%%%%%%%%%%%%%%%%%%%%%%%%%%%%%%%%%%%

\bigskip

\begin{Theorem}\label{mainresrednCM}
 Let $(\RR,\m)$ be a Noetherian local ring of dimension $d \ge 1$ with infinite residue field. For an $\m$-primary ideal $I$ and a minimal reduction $J$ of $I$, there exists a minimal reduction $Q$ of $I$ such that
 \[ \mbox{\rm red}_{Q}(I) \leq \max \{d {\cdot} \l(\RR/J) - 2d+1, 0\}.\]
\end{Theorem}

\demo Let us start with a minimal reduction  $J=(x_1, \ldots, x_{d})$
  of $I$. Let $L= (x_1, \ldots, x_{d-1})$. Then \[ \nu(I^n) \leq \nu(L^n)
+ \nu(I^n/L^n).  \] We need to estimate ${\ds \nu(I^n/L^n)}$.
%%%%%%%%%%%%%%%%%%%%%%%%%%%%%%%%%%%%%%%%%%%%%%%%%%%%%%%%%%%%%%%%%%%%%%%%%%%%%%%%
% where $I^n/L^n$ is an ideal of the $1$-dimensional local ring $\RR/L^n$.
%%%%%%%%%%%%%%%%%%%%%%%%%%%%%%%%%%%%%%%%%%%%%%%%%%%%%%%%%%%%%%%%%%%%%%%%%%%%%%%%%%
Set
\[ M_{i} = \frac{(I^n + L^{n-i+1}) \cap L^{n-i}}{L^{n-i+1}}  \quad \mbox{\rm and} \quad
N_{i}=\frac{I^n + L^{n-i+1}}{ L^{n-i+1}}. \]
Then  we obtain the following series of exact sequences :
\[ 0 \lar M_{i} \lar N_{i} \lar N_{i+1} \lar 0,  \]
where $i=1, \ldots, n-1$. Note that for each $i$, $M_{i}$ is a submodule of ${\ds L^{n-i}/ L^{n-i+1}}$ as an $\RR/L$--module.
Hence by Lemma~\ref{lengthdim1}, for each $i=1, \ldots, n-1$,
\[ \nu(M_{i}) \leq \l(\RR/J) {\cdot}\nu(L^{n-i}/L^{n-i+1}) = \l(\RR/J){\cdot} {\ds {{d+n-2-i}\choose{d-2}}}.    \]
Since ${\ds N_{n} = \frac{I^n + L}{L} }$ is a submodule of $\RR/L$, by Lemma~\ref{lengthdim1}, we get
\[ \nu(N_{n}) \leq \l(\RR/J).    \]
Therefore we obtain
\[ \begin{array}{lll}
\nu(I^n) &\leq & \nu(L^n) + \nu(I^n/L^n) \\ && \\
         % &\leq & \nu(L^n) + \nu(M_{1}) + \nu(N_{2}) \\ && \\
         %&\leq & \nu(L^n) + \nu(M_{1}) + \nu(M_{2}) +\nu(N_{2}) \\ && \\
         % &     & \vdots \\ && \\
         &\leq & \nu(L^n) + \nu(M_{1}) + \nu(M_{2}) + {\cdots} + \nu(M_{n-1})+ \nu(N_{n}) \\ && \\
         &\leq & {\ds {{d+n-2}\choose{d-2}} + \l(\RR/J) {\cdot}\sum_{i=1}^{n}  {{d+n-2-i}\choose{d-2}}} \\ && \\
         &=& {\ds {{d+n-2}\choose{d-2}} + \l(\RR/J){\cdot} {{d+n-2}\choose{d-1}}.   }
\end{array}
\]
Recall that if
\[ {\ds  \nu(I^n) < {{n+d}\choose{d}},  } \]
then there is a minimal reduction $Q$ of $I$ such that $\red_{Q}(I) \leq n-1$ (\cite{ES}, \cite[Theorem 2.36]{icbook}).
Hence by solving the inequality
\[{\ds   {{d+n-2}\choose{d-2}} + \l(\RR/J) {{d+n-2}\choose{d-1}}   <  {{n+d}\choose{d}},  }  \]
we obtain the desired relation.
\QED

\bigskip

%%%%%%%%%%%%%%%%%%%%%%%%%%%%%%%%%%%%%%%%%%%%%%%%%%%%%%%%%%%%%%%%%%%%%%%%%%%%%%%%%%%%%%%%%%%%%%%%%%%%%%
% \noindent Proof of the last equality by induction on $n$. We want to show that
% \[{\ds \sum_{i=1}^{n}  {{d+n-2-i}\choose{d-2}} =  {{d+n-2}\choose{d-1}}.     }    \]
% It is clear when $n=1$. Suppose that $n \geq 2$. Then
% \[\begin{array}{lll}
% {\ds \sum_{i=1}^{n+1}  {{d+n-1-i}\choose{d-2}}} &= & {\ds {{d+n-2}\choose{d-2}} + {{d+n-3}\choose{d-2}} + {\cdot}s + {{d-2}\choose{d-2}}} \\ && \\ &=& {\ds
% {{d+n-2}\choose{d-2}} +  \sum_{i=1}^{n}  {{d+n-2-i}\choose{d-2}}  } \\ && \\
% &=& {\ds {{d+n-2}\choose{d-2}} +  {{d+n-2}\choose{d-1}}  } \\ && \\
% &=& {\ds {{d+n-1}\choose{d-1}}  }
% \end{array}
% \]
%%%%%%%%%%%%%%%%%%%%%%%%%%%%%%%%%%%%%%%%%%%%%%%%%%%%%%%%%%%%%%%%%%%%%%%%%%%%%%%%%%%%%%%%%%%%%%%%%

\bigskip

\begin{Corollary}
Let $(\RR, \m)$ be a Noetherian local ring of dimension $d \ge 1$ and infinite residue field. Let $Q$ be a minimal reduction of $\m$ such that $\mbox{\rm red}_{Q}(\m) = \mbox{\rm red}(\m)$.
Then
\[ e_{1}(\m) \leq e_{1}(\m) - e_{1}(Q) \leq \l(\RR/(Q:\m) ) {\cdot} \left[{\ds {{ \nu(\m)+\l(\RR/Q)d -3d+1}\choose{\nu(\m)-d}}} -1 \right]. \]
\end{Corollary}

\begin{Remark}{\rm It is worthwhile to point out that there are other
known bounds for the reduction number of an ideal in terms of some of
its Hilbert coefficients. One of these is a bound proved
by M. E. Rossi (\cite[Corollary 1.5]{Rossi00}): If $(\RR, \m)$ is a Cohen-Macaulay
local ring
of dimension at most $2$ then for any $\m$--primary ideal $I $ with a minimal
reduction $Q$
\[ \mbox{\rm red}_Q(I)\leq e_1(I) - e_0(I) + \l(\RR/I) + 1.
\] Several open questions arise. Foremost whether it extends to
higher dimensional Cohen-Macaulay rings (with a correction term
depending on the dimension). Another question is which offsetting terms should be added
in the non Cohen-Macaulay case. For instance, in dimension $2$
whether the addition of $-e_1(Q)$, a term that can be considered a
non Cohen--Macaulayness penalty, would give a valid bound. 
}
\end{Remark}

\section{Normalization}

The following observation shows how the special fiber of the
normalization impacts $e_0(I)$. Of course,  more interesting issue
would be to obtain relationships going the other way.

%\begin{Conjecture}{\rm Let $(\RR,\m)$ be an unmixed local domain and
%$I$ an $\m$-primary ideal. Then $\overline{e}_1(I)\geq 0$. }
%\end{Conjecture}

%The motivation is partly based on the fact that the integral closure
%of integral domains--at least the geometric ones--is also the
%integral closure of a Cohen-Macaulay ring. (Maybe this is silly...)
%In any event one should try to understand $\overline{e}_1$.

\begin{Proposition}
Let $(\RR,\m)$ be a normal local domain and let $I$ be an $\m$-primary
ideal. Suppose that $\cl{\rees} = \bigoplus_{n=0}^{\infty} C_{n}$ is finite
over $\rees =\rees (I)$. We denote by $\cl{f}_0(I)$ the multiplicity of  $\cl{\rees}/\m \cl{\rees}$.
Then  \[e_0(I) \leq \min\{f_0(I){\cdot} \l(\RR/I),\;\; \cl{f}_0(I){\cdot}
\lambda(\RR/\cl{I})\}.\]
\end{Proposition}

\demo We first observe that $C_{n+1}= IC_n= \cl{I}C_n$, for $n \gg 0$.
In particular, in that range, $C_{n+1}\subset \m C_n$.  Consider now
the corresponding exact sequence \[ 0 \rar \m C_n/C_{n+1} \lar
C_n/C_{n+1} \lar C_n/\m C_n \rar 0. \] Counting multiplicities, we
have \[e_0(I) \leq \deg(\m \cl{\rees}/\cl{I}
\; \cl{\rees}) + \deg(\cl{\rees}/\m \cl{\rees})\leq \cl{f}_0(I)(
\l(\m/\cl{I})+1) = \cl{f}_0(I){\cdot} \lambda(\RR/\cl{I})\] as desired. The other inequality, $e_0(I)\leq
f_0(I){\cdot} \l(\RR/I)$, has a similar proof.  \QED

\end{document}